\documentclass[12pt]{article}

\usepackage{fullpage}
\usepackage{xcolor}
\usepackage{amsmath}
\usepackage{amssymb}
\usepackage{amsthm,mathtools}
\usepackage{bbm,bm}
\usepackage{hyperref}
\usepackage{theoremref}
\usepackage{enumitem}
\usepackage{stmaryrd} 
\usepackage[colors]{optsys}

\def\reals{\mathbb{R}}
\def\ereals{\overline{\reals}}

\def\comp{\mathop{\text{\scriptsize\raise 1pt \hbox{$\circ$}}}}
\def\infconv{\mathop{\text{\scriptsize\raise 1pt \hbox{$\square$}}}}

\def\argmin{\mathop{\rm argmin}\limits}
\def\argmax{\mathop{\rm argmax}\limits}

\def\minimize{\mathop{\rm minimize}\limits}
\def\maximize{\mathop{\rm maximize}\limits}

\def\st{\mathop{\rm subject\ to}}

\def\dom{\mathop{\rm dom}\nolimits}

\def\ovr{\mathop{\rm over}\ }

\def\Tr{\mathop{\rm Tr}}

\def\upto{{\raise 1pt \hbox{$\scriptstyle \,\nearrow\,$}}}
\def\downto{{\raise 1pt \hbox{$\scriptstyle \,\searrow\,$}}}

\def\epi{\mathop{\rm epi}\nolimits}

\def\FF{(\F_t)_{t=0}^T}

\def\F{{\cal F}}

\def\N{{\cal N}}

\def\S{{\cal S}}

\newcommand{\utility}{\upsilon}

\newcommand{\nb}[3]{
  {\colorbox{#2}{\bfseries\sffamily\tiny\textcolor{white}{#1}}}
  {\textcolor{#2}{$\blacktriangleright${#3}$\blacktriangleleft$}}}

\newcommand{\mdl}[1]{\nb{Michel}{blue}{#1}}

\title{Optimal Operation and Valuation\\ of Electricity Storages \footnote{This
    research benefited from the support of the FMJH Program Gaspard Monge for
    optimization and operations research and their interactions with data
    science.}}

\author{Jean-Philippe Chancelier$^\dagger$ \and
  Michel De Lara\thanks{CERMICS, \'Ecole nationale des ponts et chaussées, IP Paris, France, {jean-philippe.chancelier, michel.delara}@enpc.fr}
  \and
  Fran\c{c}ois Pacaud\thanks{Centre Automatique et Systèmes, Mines Paris-PSL, 60 boulevard Saint Michel, 75006 Paris, francois.pacaud@mineparis.psl.eu}
  \and
  Teemu Pennanen\thanks{Department of Mathematics, King's College London, Strand, London, WC2R 2LS, United Kingdom, teemu.pennanen@kcl.ac.uk}
  \and
  Ari-Pekka Perkki\"o\thanks{Mathematics Institute, Ludwig-Maximilian University of Munich, Theresienstr. 39, 80333 Munich, Germany, a.perkkioe@lmu.de}}

\begin{document}
\maketitle

\begin{abstract}
  This paper applies computational techniques of convex stochastic optimization
  to optimal operation and valuation of electricity storages in the face of
  uncertain electricity prices. Our approach is applicable to various
  specifications of storages, and it allows for e.g.\ hard constraints on
  storage capacity and charging speed. Our valuations are based on the
  indifference pricing principle, which builds on optimal trading strategies and
  calibrates to the user's initial position, market views and risk
  preferences. We illustrate the effects of storage capacity and charging speed
  by numerically computing the valuations using stochastic dual dynamic
  programming.
\end{abstract}

\section{Introduction}

The growing proportion of renewable energy generation increases the
uncertainties and seasonalities of supply and the price of electricity. This
creates an incentive to store energy~\cite{weitzel2018energy} for the purposes
of dynamic trading.  This article proposes mathematical models and computational
tools for optimal operation and valuation of energy storages in the face of
uncertain electricity prices and investment returns. The presented techniques
apply to various physical specifications of storage and they yield
market-consistent valuations of storage facilities. Economically sensible
valuation and optimization techniques will facilitate the expansion of the
overall storage capacity which will be necessary for the development of
resilient power systems based on renewables~\cite{yu2017stochastic}.

We first develop a stochastic optimization model for finding optimal strategies
for charging and discharging a given energy storage subject to physical
constraints. The aim is to optimize a user-specified risk measure of the
terminal wealth which depends on the underlying risk factors and the
decision-maker's strategy. Realistic models of the storage result in stochastic
optimization problems that fail to be analytically solvable. This is partly due
to the fact that a physical storage has capacity constraints which amount to
state constraints in the corresponding optimal control formulation. We will thus
employ computational techniques for the numerical solution of the resulting
optimization problems.

Being able to solve the optimal storage management problem numerically allows us
to compute {\em indifference prices} for the storages. Indifference pricing has
been extensively studied in pricing financial products in incomplete markets;
see e.g.\ the collection~\cite{car9} and the references there. In~\cite{ptw9},
the authors applied indifference pricing to price a production facility that is
controlled by switching between a finite number of states. Closer to the present
work is~\cite{ccgv}, where utility indifference pricing was applied to various
energy contracts with nontrivial optionalities. A particularly interesting
example is a ``virtual storage" contract, which resembles physical storages in
terms of the generated cashflows. The numerical computations in~\cite{ccgv} are
based on finite difference approximations of the associated
Hamilton-Jacobi-Bellman (HJB) partial differential equation. As is typical of
HJB-based approaches in optimal control, such techniques do not allow for state
constraints so they are not appropriate to physical storages.
The recent paper~\cite{pv21} studies optimal operation of a pumped hydroelectric
storage with state constraints. The authors assume positive electricity prices
and linear utility functions. They treat state constraints by the ``level set
approach" where nondifferentiable penalty functions are added to the
objective. Similar techniques were applied by~\cite{gpw23} in the context of
mean field games with some applications to renewable energy.

The optimal management of power systems under uncertainty has been recently
surveyed in~\cite{roald2023power}, with the stochastic dual dynamic programming
(SDDP) algorithm identified as one of the most promising method for the
management of energy
storage~\cite{lohndorf2013optimizing,megel2015stochastic,lw21,zheng2022arbitraging}.
On the contrary to model predictive control (MPC)-based
approach~\cite{arnold2011model}, SDDP takes explicitly into account the future
uncertainties when optimizing decisions.  Compared to other classical approaches
in stochastic optimal control~\cite{ccgv}, SDDP has the advantage of allowing
for state constraints, necessary for realistic descriptions of storages which
usually involve hard constraints on the amount of stored energy. Unlike
\cite{lw21}, our model includes the possibility of lending or borrowing money in
a cash account/money market. This allows us to avoid ``discounting" future
cashflows and it leads to a more natural notion of an indifference price of a
storage in terms of money paid at the time of purchase. Moreover, thanks to the
explicit modeling of the agent's cash position, we avoid the delicate
``asset-backed trading strategies" of~\cite{lw21} that, at a given point~$t$ in
time, involve planned future decisions in addition to the decisions to be
implemented at time~$t$.

This paper builds on the stochastic dual dynamic programming (SDDP) algorithm to
find approximately optimal storage management strategies over multiple
periods. The presented approach is quite general and it can be applied to
various storage and market specifications. The approach is illustrated on a
fairly simple problem instance over a 30-day period. We provide numerical
illustrations of storage valuations and show, in particular, how the physical
characteristics of the storage affect the valuations. The valuations react
consistently also to variations in the user's views concerning the future
development of electricity prices. The presented approach allows also for
extensions where the agent might have existing storage, production facilities or
consumption.

The rest of this paper is organized as follows.  In Sect.~\ref{sec:model}, we
present a generic mathematical formulation of the optimal management of a
storage operating on an energy market, and how to deduce an indifference pricing
(valuation).
In Sect.~\ref{sec:oc}, we outline how to adapt the Stochastic dual dynamic
programming (SDDP) algorithm to the mathematical formulation above, but in the
so-called Markovian noise case.
Sect.~\ref{sec:example} describes a simple instance of the general storage
management model, followed by numerical results in Sect.~\ref{sec:comp}.

\section{Optimal storage management and valuation}
\label{sec:model}

In Section~\ref{Optimal_storage_management}, we present a mathematical
formulation of the problem of optimal management of an electricity storage
connected to an energy market.
In Section~\ref{sec:idp}, we recall the notion of indifference price and how it
can be deduced from the optimization problem above.

\subsection{Optimal storage management}
\label{Optimal_storage_management}

Consider an agent who has an electricity storage and can trade in the
electricity market over a given discrete time span~$\ic{0,T}=\{0,1,\ldots,T-1,T\}$,
where $T\geq 1$ is a natural number (horizon). The agent aims to use the
storage optimally to maximize his terminal wealth obtained by trading
electricity. Any money received by selling energy is invested in financial
markets, and any energy bought will be financed by liquidating financial
investments. Finding the optimal way to operate the storage is nontrivial due to
the uncertainty concerning the prevailing electricity prices and investment
returns in the financial market. We will model the electricity prices and
investment returns as stochastic processes, and formulate the problem as a
stochastic optimization problem. Over shorter planning horizons, the
uncertainties in investment returns can often be ignored, but the uncertainty of
electricity prices may still be significant.

We will formulate the storage management problem as a general {\em convex
  stochastic optimal control} problem over finite discrete time and a general
probability space $(\Omega,\F,P)$ with a filtration $\FF$ (an increasing sequence of
sub-$\sigma$-algebras of $\F$).  We denote by $\EE[\cdot]$ the mathematical expectation
under probability~$P$.  We will allow for the general control format
of~\cite[Section~6.2]{ccdmt} where, at each stage $t$, the system may be
controlled by both $\F_{t-1}$-measurable and $\F_t$-measurable controls; see
also~\cite{pp25a}. Such formulations arise naturally, e.g.\ when trading both in
the day-ahead auctions as well as intraday markets. Other instances can be found
in~\cite{pdccjr}.

Consider the optimal control problem
\begin{subequations}
  \label{oc}
  \begin{align}
    &\minimize & \EE\Bigg[\sum_{t=0}^{T-1}&  L_t(X_t,U^b_t,U^a_{t+1})+ J(X_T)\Bigg]\ \ovr\ (X,U^b,U^a)\in\N,\\ 
    &\st\ & X_0 &= x_0,\\
    & & X_{t+1} &= F_{t}(X_t,U^b_t,U^a_{t+1})\quad t=0,\dots,T-1\ a.s.,
                  \label{oc_c}
\end{align}
where the state $X_t$ and the controls\footnote{The superscript~$b$ in~$U^b_t$ stands for ``before'' (in the sense of a decision taken at the beginning of the time interval~$[t,t+1[$), whereas the superscript~$a$ in~$U^a_t$ stands for ``after'' (in the sense of a decision taken at the end of the time interval~$[t,t+1[$).} $U^b_t$ and $U^a_t$ are random variables taking values in $\reals^N$, $\reals^{M^b}$ and $\reals^{M^a}$. The functions $L_t:\reals^N\times\reals^{M^b}\times\reals^{M^a}\times\Omega\to\ereals$ and $J:\reals^N\times\Omega\to\ereals$ are extended real-valued convex normal integrands\footnote{Recall that a function $f:\reals^n\times\Omega\to\ereals$ is a {\em normal integrand} if the epigraphical mapping $\omega\mapsto\epi f(\cdot,\omega)$ is a closed-valued and measurable.} and, for each $t$,
\begin{equation*}
  \label{eq:system_equations}
  F_{t}(x_t,u^b_t,u^a_{t+1},\omega) := A_{t+1}(\omega)x_t+B^b_{t+1}(\omega)u^b_t+B^a_{t+1}(\omega)u^a_{t+1}+W_{t+1}(\omega),
\end{equation*}
\end{subequations}
where $A_{t+1}$, $B^b_{t+1}$ and $B^a_{t+1}$ are $\F_{t+1}$-measurable random matrices of appropriate dimensions and $W_{t+1}$ are $\F_{t+1}$-measurable random vectors. In \eqref{oc}, we have omitted the $\omega$ from the notation for brevity. The set $\N$ denotes the space of adapted stochastic processes, i.e.\ those where $(X_t,U^b_t,U^a_t)$ is $\F_t$-measurable for each $t=0,\ldots,T$. The vector $x_0\in\reals^N$ is a given initial state. Accordingly, we will assume that $\F_0$ only consists of sets of measure zero or one, so that any $\F_0$-measurable function is almost surely constant.

Problem~\eqref{oc} is quite general and can accommodate various kinds of
storages and market specifications. Besides storage management and electricity
trading, it can include production facilities and/or consumption of
electricity. The computational techniques discussed below apply to such
situations as well. Both production and consumption may create stronger
incentives to buy storages. The indifference pricing approach described in
Section~\ref{sec:idp} below is well suited to account for such
modifications.



\subsection{Indifference pricing}\label{sec:idp}

Indifference pricing is a general principle for defining economically meaningful
values of products in incomplete markets where prices cannot be determined by
arbitrage arguments alone. We define the {\em indifference price} for
buying/renting a storage as the maximum amount an agent could pay for it without
worsening his well-being as measured by the optimum value of
problem~\eqref{oc}. In order to emphasize its dependence on the initial state
$x_0$, we will denote the optimum value of problem~\eqref{oc} by $\varphi(x_0)$.

For an agent who does not have access to a storage, problem~\eqref{oc} becomes
\begin{equation}\label{oc0}
  \begin{aligned}
    &\minimize & \EE\Bigg[\sum_{t=0}^{T-1}&  L_t(X_t,U^b_t,U^a_{t+1})+ J(X_T)\Bigg]\ \ovr\ (X,U^b,U^a)\in\N,\\ 
    &\st\ & X_0 &= x_0,\\
    & & X_{t+1} &= F_{t}(X_t,U^b_t,U^a_{t+1})\quad t=0,\dots,T-1\ a.s.,\\
    & &  (X_t&, U^b_t,U^a_{t+1}) \in D_t,\quad t=0,\dots,T-1\ a.s.,
  \end{aligned}
\end{equation}
where $D_t$ describes a constraint that prohibits using the storage at time
$t$. More concretely, if the state variable can be decomposed as
$X=(X^e,\tilde X)$, where $X^e$ denotes the amount of stored energy and
$\tilde X$ denotes the other state variables in the storage management
problem~\eqref{oc}, then one could simply set
$D_t=\{0\}\times\reals^{N-1}\times\reals^{M^b}\times\reals^{M^a}$, which says that the level of
the storage has to remain at zero; see Section~\ref{sec:tom} for a simple
example.

We denote the optimum value of~\eqref{oc0} by $\psi(x_0)$. The {\em indifference
  price} $\pi$ of the storage is then defined as
\[
  \pi := \sup\bset{\alpha\in\reals}{\varphi(x_0-\alpha p)\le\psi(x_0)}
  \eqfinv
\]
where $p\in\reals^N$ is the unit vector that has $1$ in the cash component and
zero elsewhere.

\begin{proposition}\label{prop:idp}
  The optimum value $\varphi(x_0)$ of~\eqref{oc} is convex as a function
  of~$x_0$. If the function $\alpha\mapsto\varphi(x_0-\alpha p)$ is strictly decreasing, then the
  indifference price $\pi$ is the unique solution to the equation
  \begin{equation}\label{eq:root}
    \varphi(x_0-\alpha p) = \psi(x_0).
  \end{equation}
\end{proposition}

\begin{proof}
  The convexity follows from the fact that the essential objective of~\eqref{oc}
  is jointly convex in $x_0$ and $(X,U^b,U^a)$; see e.g.\ Theorem~1 and
  Example~1 in~\cite{roc74}. The convexity of $\varphi$ implies that it is continuous
  with respect to the interior of its domain
  \[
    \dom\varphi:=\{x_0\in\reals\,|\, \varphi(x_0)<+\infty\};
  \]
  see~\cite[Theorem~10.1]{roc70a}. The second claim thus follows from the mean
  value theorem.
\end{proof}

Finding a solution of the equation~\eqref{eq:root} is a nontrivial problem in
general since a single evaluation of the function $\varphi$ requires the numerical
solution problem~\eqref{oc}. The task is facilitated, however, by the fact that
$\varphi$ is nondecreasing and convex. In some situations, a solution can be found
analytically after two solves of problem~\eqref{oc}; see Section~\ref{sec:tom}.

\section{Stochastic dual dynamic programming}
\label{sec:oc}

In general, problem~\eqref{oc} does not allow for analytical solutions, so we
will resort to computational techniques based on dynamic programming, like
exposed in Section~\ref{Stochastic_dynamic_programming}.  Then, we outline the
so-called Markovian case, together with discretization of Markov processes, in
Section~\ref{subsec:dmp}.  Finally, we describe the stochastic dual dynamic
programming algorithm in Section~\ref{subsec:sddp}.

\subsection{Stochastic dynamic programming}
\label{Stochastic_dynamic_programming}

We will assume throughout that there is a stochastic process
$\xi=(\xi_t)_{t=0}^T$ such that $L_{t-1}$, $ A_t$, $ B^b_t, B^a_t$ and $ W_t$
in~\eqref{eq:system_equations} depend on $\omega$ only through
$\xi^t:=(\xi_s)_{s=0}^t$ and $J_T$ depends on $\omega$ only through $\xi^T$. We will also
assume that $\xi_t$ takes values in $\reals^{d_t}$ and we denote
\[
d^t:=d_0+\cdots+d_t
\]
and $d=d^T$. Without loss of generality, we will assume that $\Omega=\reals^d$ and
that $\xi_t(\omega)$ is the projection of $\omega$ to its $t$th component and that $\F_t$ is the $\sigma$-algebra on $\reals^d$ generated by the projection
$\xi\mapsto\xi^t:=(\xi_0,\ldots,\xi_t)$. In particular, $P$ is the distribution of~$\xi$. With a slight misuse of notation we will write 
$L_t(x_t,u^b_t,u^a_{t+1},\xi^{t+1})$ instead of $L_t(x_t,u^b_t,u^a_{t+1},\omega)$ and similarly for other $\F_{t+1}$-measurable random elements. 

It follows that the $\F_t$-conditional expectation of an $\F_{t+1}$-measurable quasi-integrable random variable $\xi\mapsto\eta(\xi^{t+1})$ has a pointwise representation
\[
\EE\big[\eta(\xi^{t+1})| \xi^t\big] := \int_{\reals^{d_{t+1}}}\eta(\xi^t,\xi_{t+1})d\nu_t(\xi^t,d\xi_{t+1})\quad\forall \xi^t\in\reals^{d^t},
\]
where $\nu_t$ is a {\em regular $\xi^t$-conditional distribution of $\xi_{t+1}$}; see, e.g.,~\cite[Section~2.1.5]{pp24}. Note that $\xi^{t+1}=(\xi^t,\xi_{t+1})$. We say that a function $h:\reals^n\times\reals^{d^{t+1}}\to\ereals$ is {\em lower bounded} if there is an integrable function $m:\reals^{d^{t+1}}\to\reals$ such
that $h(x,\xi^{t+1})\ge m(\xi^{t+1})$ for all
$(x,\xi^{t+1})\in\reals^n\times\reals^{d^{t+1}}$. Given a lower bounded measurable
function $h:\reals^n\times\reals^{d^{t+1}}\to\ereals$, we denote
\[
\EE[h(x,\xi^{t+1})|\xi^t] := \int_{\Xi_{t+1}}h(x,\xi^t,\xi_{t+1})d\nu_t(\xi^t,d\xi_{t+1})\quad\forall (x,\xi^t)\in\reals^n\times\reals^{d^t}.
\]
  
We say that a sequence of lower bounded measurable functions
$J_t:\reals^N\times\reals^{d^t}\to\ereals$ solves the Bellman equations for~\eqref{oc}
if $J_T=J$ and 
\begin{equation}\label{be}
  J_t(x_t,\xi^t) = \inf_{u^b_t}\left\{\EE\left[\inf_{u^a_{t+1}}\{L_t(x_t,u^b_t,u^a_{t+1},\xi^{t+1}) + J_{t+1}(F_{t}(x_t,u^b_t,u^a_{t+1},\xi^{t+1}),\xi^{t+1})\}\Bigg|\xi^t\right]\right\}
\end{equation}  
for all $t=0,\ldots,T-1$ and $P$-almost every $\xi$.
Bellman equations of the above form~\eqref{be}, together with optimal feedbacks, can be derived\footnote{%
  In \cite[Sect.~6]{ccdmt}, Proposition~13 can be extended to the case where instantaneous costs, dynamics
  and conditional distribution of the next uncertainty all three depend upon the whole sequence of past uncertainties (and not only upon the last uncertainty).
  The proof requires an adaptation of Appendix~D in the preprint version of~\cite{ccdmt} (which contains additional
material compared to the published version~\cite{ccdmt}). 
We can then obtain a variant of \cite[Equation~(42)]{ccdmt} that corresponds to Equation~\eqref{be} above.
Then, optimal feedbacks can be deduced from \cite[Subsect.~3.3]{ccdmt}.
}
from the framework of~\cite[Sect.~6]{ccdmt} (and \cite[Subsect.~3.3]{ccdmt} for optimal feedbacks), 
which first analyzed 
Bellman equations in the case of nonlinear system dynamics
with decision-hazard-decision information structures.
The equations are well-motivated by the following optimality conditions from \cite[Theorem~1]{pp25a}.

\begin{theorem}[Optimality principle]\thlabel{opoc}
If $(J_t)_{t=0}^T$ is a solution of~\eqref{be} such that each
$J_t:\reals^N\times\Xi^t\to\ereals$ is a lower bounded convex normal integrand, then the optimum value of~\eqref{oc} equals $J_0(x_0)$ and an $(\bar X,\bar U^b,\bar U^a)\in\N$ solves~\eqref{oc} if and only if $\bar X_0=x_0$ and
\begin{align*}
\bar U^b_t(\xi^t)&\in \argmin_{u^b_t\in\reals^{M^b}}
      \left\{\EE\left[\inf_{u^a_{t+1}\in\reals^{M^a}}\{L_t(\bar X_t,u^b_t,u^a_{t+1},\xi^{t+1}) + J_{t+1}(F_{t}(\bar X_t,u^b_t,u^a_{t+1},\xi^{t+1}),\xi^{t+1})\}\Bigg|\xi^t\right]
      \right\}\eqfinv
    \\
    \bar U^a_{t+1}(\xi^{t+1})
    &\in \argmin_{u^a_{t+1}\in\reals^{M^a}}\{L_t(\bar X_t(\xi^t),\bar U^b_t(\xi^t),u^a_{t+1},\xi^{t+1})
     + J_{t+1}(F_{t}(\bar X_t(\xi^t),\bar U^b_t(\xi^t),u^a_{t+1},\xi^{t+1}),\xi^{t+1})\},
    \\
    \bar X_{t+1}(\xi^{t+1})
    &= F_{t}(\bar X_t(\xi^t),\bar U^b_t(\xi^t),\bar U^a_{t+1}(\xi^{t+1}),\xi^{t+1})
  \end{align*}
  for all $t=0,\ldots,T-1$ and $P$-almost every $\xi$.
\end{theorem}

The following gives sufficient conditions for the existence of solutions to the
Bellman equations~\eqref{be}. Given a normal integrand
$h:\reals^n\times\Omega\to\ereals$, the function~$h^\infty$, obtained from $h$ by defining
$h^\infty(\cdot,\omega)$, for each $\omega\in\Omega$, as the {\em recession function} of
$h(\cdot,\omega)$, is a positively homogeneous convex normal integrand; see
e.g.~\cite[Example~14.54a]{pp24} or~\cite[Theorem~1.39]{pp24}. The following
is~\cite[Theorem~2]{pp25a} with the first stage running cost redefined as
$+\infty$ unless $X_0=x_0$.

\begin{theorem}[Existence of solutions]\label{existoc}
  Assume that the set
  \begin{multline*}
    \Big\{(X,U^b,U^a)\in\N \, \Big\vert \,
    \sum_{t=0}^{T-1}  L^\infty_t(X_t,U^b_t,U^a_{t+1})+ J^\infty(X_T)\le 0,\\
    X_0=0\eqsepv X_t= A_tX_{t-1} +  B^b_tU^b_{t-1} +  B^a_tU^a_t\eqsepv t=1,\dots,T\ a.s.
    \Big\}
  \end{multline*}
  is linear. Then the control problem~\eqref{oc} has an optimal solution and the
  associated Bellman equations~\eqref{be} have a unique solution $(J_t)_{t=0}^T$
  of lower bounded convex normal integrands.
\end{theorem}

The following is from~\cite[Theorem~3]{pp25a}.

\begin{theorem}[Markovianity]\label{thm:markov}
  Let $(J_t)_{t=0}^T$ be a sequence of lower bounded normal integrands that
  solve~\eqref{be}. Assume that $\xi=(\xi_t)_{t=0}^T$ is a Markov process and that
  $L_{t-1}$, $A_t$, $B^b_t$, $B^a_t$ and $W_t$ depend on $\xi$ only through
  $(\xi^t,\xi^{t-1})$. Then $J_t$ depends on $\xi^t$ only through $\xi_t$. If
  $\xi$ is a sequence of independent random variables and if $L_{t-1}$, $A_t$,
  $B^b_t$, $B^a_t$ and $W_t$ do not depend on $\xi^{t-1}$, then $J_t$ do not
  depend on $\xi$ at all.
\end{theorem}

\subsection{Discretization of Markov processes}\label{subsec:dmp}

In general, the Bellman equations~\eqref{be} do not allow for analytical
solutions. We will therefore resort to numerical approximations and assume, as
in \autoref{thm:markov}, that there is a Markov process $\xi=(\xi_t)_{t=0}^T$ such
that $L_{t-1}$, $A_t$, $B_t$ and $W_t$ depend on $\xi$ only through
$(\xi_{t-1},\xi_t)$. 
We will first approximate the Markov process~$\xi$ by a finite-state Markov chain,
and then use the Stochastic Dual Dynamic Programming (SDDP) algorithm to solve
the discretized problem; see Section~\ref{subsec:sddp} below.

We will discretize $\xi$ using integration quadratures as
in~\cite{tauchen1991quadrature} in the context of option pricing. To that end,
we assume that the $\xi_t$-conditional distribution of $\xi_{t+1}$ has density
$p_t(\cdot|\xi_t)$. Given a strictly positive probability density $\phi_{t+1}$ on
$\reals^{d_{t+1}}$, the conditional expectation of a random variable
$\xi\mapsto\eta(\xi^{t+1})$ will be approximated by
\begin{align*}
  \EE[\eta |\xi_t] &= \int_{\reals^{d_{t+1}}}\eta(\xi_{t+1}) p_t(\xi_{t+1}|\xi_t)d\xi_{t+1} \\
              &=\int_{\reals^{d_{t+1}}}\eta(\xi_{t+1})\frac{p_t(\xi_{t+1}|\xi_t)}{\phi_{t+1}(\xi_{t+1})}\phi_{t+1}(\xi_{t+1})d\xi_{t+1}\\
              &\approx \sum_{i=1}^N\eta(\xi_{t+1}^i)\frac{p_t(\xi_{t+1}^i|\xi_t)}{\phi_{t+1}(\xi_{t+1}^i)}w_{t+1}^i,
\end{align*}
where $(\xi_{t+1}^i,w_{t+1}^i)_{i=1}^N$ is a quadrature approximation of the
density $\phi_t$. A natural choice for $\phi_{t+1}$ would be the unconditional density
of $\xi_{t+1}$ provided it is available in closed form. Constructing such
approximations for $t=0,\ldots,T$ and conditioning each $\xi_{t+1}$ with the quadrature
points $\xi_t^i$ obtained at time $t$ results in a finite state Markov chain
where, at each~$t$, the process~$(\xi_t)_{t=0}^T$ takes values in the finite set
$\{\xi_t^j\}_{j=1}^N$ and the transition probability from $\xi_t^j$ to
$\xi_{t+1}^i$ is given by
\[
p_t^{ji}:=\frac{p_t(\xi_{t+1}^i|\xi^j_t)}{\phi_{t+1}(\xi_{t+1}^i)}w_{t+1}^i.
\]
Sufficient conditions for the weak convergence of such discretizations are discussed in~\cite{kar75}.

\subsection{Stochastic dual dynamic programming algorithm}\label{subsec:sddp}

Stochastic dual dynamic programming is a computational technique for the
solution of the Bellman equations in convex stochastic optimization problems
governed by finite-state Markov chains. The SDDP algorithm was introduced
in~\cite{Pereira-Pinto:1991} for stagewise independent processes $\xi$ and later
extended to finite-state Markov processes $\xi$ in~\cite{im96}; see
also~\cite{pm12,lohndorf2019modeling} and the references there.

The SDDP algorithm is based on successive approximation of the cost-to-go functions $(J_t)_{t=0}^T$ that solve~\eqref{be} by polyhedral lower approximations. For control problems in the format of~\eqref{oc}, the algorithm
proceeds as follows\footnote{Note that, when the underlying process $\xi$ is a finite-state Markov chain, the conditional expectations below are finite sums weighted by the transition probabilities.}:
{\footnotesize
  \begin{enumerate}
  \item[0.]
    {\bf Initialization:} Choose convex (polyhedral) lower-approximations $J_t^0$ of the cost-to-go functions $J_t$ and set $k=0$.
  \item
    {\bf Forward pass:} Sample a path $\xi^k$ of $\xi$ and define $X^k$ by $X_0^k=x_0$ and
    \begin{align*}
U_t^{bk}&\in\argmin_{u^b_t\in\reals^{M^b}}\EE\left[\inf_{u^a_{t+1}\in\reals^{M^a}}\{L_t(X^k_t,u^b_t,u^a_{t+1}) + J^k_{t+1}(F_{t}(X^k_t,u^b_t,u^a_{t+1}))\}
                \Bigg{|}\xi^k_t\right],\\
      U_{t+1}^{ak} &\in\argmin_{u^a_{t+1}\in\reals^{M^a}}\{L_t(X^k_t,U^{bk}_t,u^a_{t+1},\xi^k_t,\xi^k_{t+1}) + J^k_{t+1}(F_{t}(X^k_t,U^{bk}_t,u^a_{t+1},\xi^k_t,\xi^k_{t+1}),\xi^k_{t+1})\},\\
      X^k_{t+1} &= F_{t}(X^k_t,U^{bk}_t,U^{ak}_{t+1},\xi^k_t,\xi^k_{t+1}).
    \end{align*}
  \item
    {\bf Backward pass:} Let $J_T^{k+1}:=J$ and, for $t=T-1,\ldots,0$,
    \begin{align*}
      \tilde J_t^{k+1}(X_t^k,\xi_t^k)
      &:= \inf_{u^b_t\in\reals^{M^b}}\EE\left[\inf_{u^a_{t+1}\in\reals^{M^a}}\{L_t(X^k_t,u^b_t,u^a_{t+1})
        + J^{k+1}_{t+1}(F_{t}(X^k_t,u^b_t,u^a_{t+1}))\}\Bigg{|}\xi^k_t\right],
      \\  
      V^{k+1}_t&\in\partial\tilde J^{k+1}_t(X_t^k,\xi_t^k)
    \end{align*}
    and let
    \[
      J_t^{k+1}(x_t,\xi_t^k) := \max\{J_t^k(x_t,\xi_t^k),\tilde J_t^{k+1}(X_t^k,\xi_t^k) + V_t^{k+1}\cdot(x_t-X^k_t)\}\quad\forall x_t\in\reals^N.
    \]
    Set $k:=k+1$ and go to 1.
  \end{enumerate}
} Above, the dot ``$\cdot$'' represents the Euclidean inner product. The minimization problems in the above algorithm are two-stage stochastic optimization problems while the original problem~\eqref{oc} has $T$ decision stages. Note also that the third minimization in the forward pass is one of the second-stage problems already solved in the second minimization. The functions $J_t^k$ generated by the algorithm are polyhedral lower approximations of the true cost-to-go functions $J_t$. It follows that the minimization
problems in the above recursions are finite-dimensional optimization problems that can be solved with standard convex optimization solvers. If the functions $L_t$ and $J$ are polyhedral, the problems are LP problems.

In the decision-observation format, where only the control $U^b$ is present, the two-stage problems reduce to a single-stage stochastic optimization problem, while in the observation-decision format without the $U^b$ control, only
deterministic single-stage problem need to be solved. This version of the algorithm is implemented e.g.\ in the Julia package described in~\cite{dowson2021sddp}.

\section{An example}\label{sec:example}

This section describes a simple instance of the general storage management
model.  Section~\ref{sec:tom} gives the problem specification and formulates it
in the general optimal control format~\eqref{oc}. Section~\ref{sec:pricemodel}
describes a simple electricity price model that captures the annual, weekly and
daily seasonalities.

The model will be used in the numerical illustrations in the following
Section~\ref{sec:comp}.  The model specification below should only be taken as
the simplest nontrivial example of the general model~\eqref{oc}.  Using the SDDP
algorithm within indifference pricing applies to any problem that can be written
in the control format~\eqref{oc}.

\subsection{The optimization model}\label{sec:tom}
We study a simple storage with a maximum capacity of $\overline X^e$ and
charging and discharging speeds of $\overline U$ and $\underline U$,
respectively. If $x$~units are stored at the beginning of period $[t,t+1[$, the
storage level will be $(1-l)x$ units at the end of the period, where the
fraction term~$l$ represents a loss per period. We describe the electricity
price uncertainty with a single stochastic process $s=(s_t)_{t+0}^T$ where $s_t$
is the unit price of electricity observed at time $t$. We assume, for
simplicity, that the same unit price applies to both purchases as well as
sales. Moreover, we assume that cash investments have a deterministic
interest~$r$ over each period $[t,t+1[$.



Assuming that the agent wants to maximize expected utility from terminal wealth, the storage management problem can be written as
\begin{equation}
  \label{sps}
  \begin{aligned}
    &\maximize\ & \besp{\utility(X^m_T)}\quad& \ovr\ (X,U)\in\N\\
    &\st & X^m_0&=x_0,\\
    & & X^e_0&=0,\\
    & & X^m_t&= (1+r) X^m_{t-1}-s_tU_t\quad t=1,\ldots,T,\\
    & & X^e_t&= (1-l) X^e_{t-1} + U_t\quad t=1,\ldots,T,\\
    & & X^e_t&\in[0,\overline X^e]\quad t=0,\ldots,T,\\
    & & U_t&\in[\underline U,\overline U]\quad t=0,\ldots,T,
  \end{aligned}
\end{equation}
where $X=(X^e,X^m)$, $X^e_t$ is the amount of stored energy, $X^m_t$ denotes
investments in money markets and $U_t$ is the amount of energy bought at
time~$t$. As usual, negative purchases are interpreted as sales. The utility
function $\utility \colon \reals\to\reals$ is assumed increasing and concave, and
$x_0$ denotes the initial cash endowment of the agent. Note that with the change
of variables $\hat X_{t+1}=X_t$, problem~\eqref{sps} is in the format
of~\eqref{oc} without the here-and-now controls $U^b$. Indeed, problem
\eqref{sps} is then an instance of \eqref{oc} with $N=2$, $M^a=1$,
\[
  J(x_T)=-\utility(x_T^m)
\]
and
\[
  L_t(x_t,u_{t+1}) = 
  \begin{cases}
    0 & \text{if $x^e_t\in[0,\overline X^e]$ and $u_{t+1}\in[\underline U,\overline U]$},\\
    +\infty & \text{otherwise}.
  \end{cases}
\]
We can thus use the numerical techniques described in Section~\ref{subsec:dmp}
and Section~\ref{subsec:sddp} above. In particular, since there are no
here-and-now controls, we can use the SDDP Julia package described
in~\cite{dowson2021sddp}.

In the numerical illustrations below, we will use the exponential utility
function
\begin{equation}\label{eq:exp}
  \utility(z) = \frac{1}{\rho}[1-\exp(-\rho z)],
\end{equation}
where $\rho>0$ is a given risk aversion parameter. Since the money market
investments are assumed to have constant interest $r$ and since the utility is
exponential, it can be shown that the indifference price given by
Proposition~\ref{prop:idp} can be expressed as
\[
  \pi = \frac{\ln(\rho\psi(x_0)+1)-\ln(\rho\varphi(x_0)+1)}{\rho(1+r)^T},
\]
where $\psi(x_0)$ is the optimum value of a problem \eqref{oc0} where use of the
storage is prohibited. Note that in the case of model \eqref{sps}, problem
\eqref{oc0} becomes
\begin{equation}
  \begin{aligned}
    &\maximize\ & &\besp{\utility(X^m_T)}\quad \ovr\ (X,U)\in\N\\
    &\st & X^m_0&=x_0,\\
    & & X^m_t&= (1+r) X^m_{t-1}\quad t=1,\ldots,T,
  \end{aligned}
\end{equation}
which has a trivial solution. Thus, the computation of the price only requires
the solution of one instance of problem~\eqref{sps}.

\subsection{A simple spot price model}\label{sec:pricemodel}
In the numerical illustrations of the following Sect~\ref{sec:comp}, we will use
a simple Markov model to describe the evolution of the electricity price
$s=(s_t)_{t=0}^T$. The model incorporates the annual, weekly and daily average
patterns observed in historical prices. We emphasize that the model developed
here is only used to illustrate the computational approach taken in
Sect~\ref{sec:comp}. The general approach is applicable to any Markov model for
the underlying risk factors. In addition to the electricity price, e.g.\ the
interest rate earned on the money market investments could be stochastic as
well.

\begin{figure}[!ht]
  \centering
  \includegraphics[width=.7\textwidth]{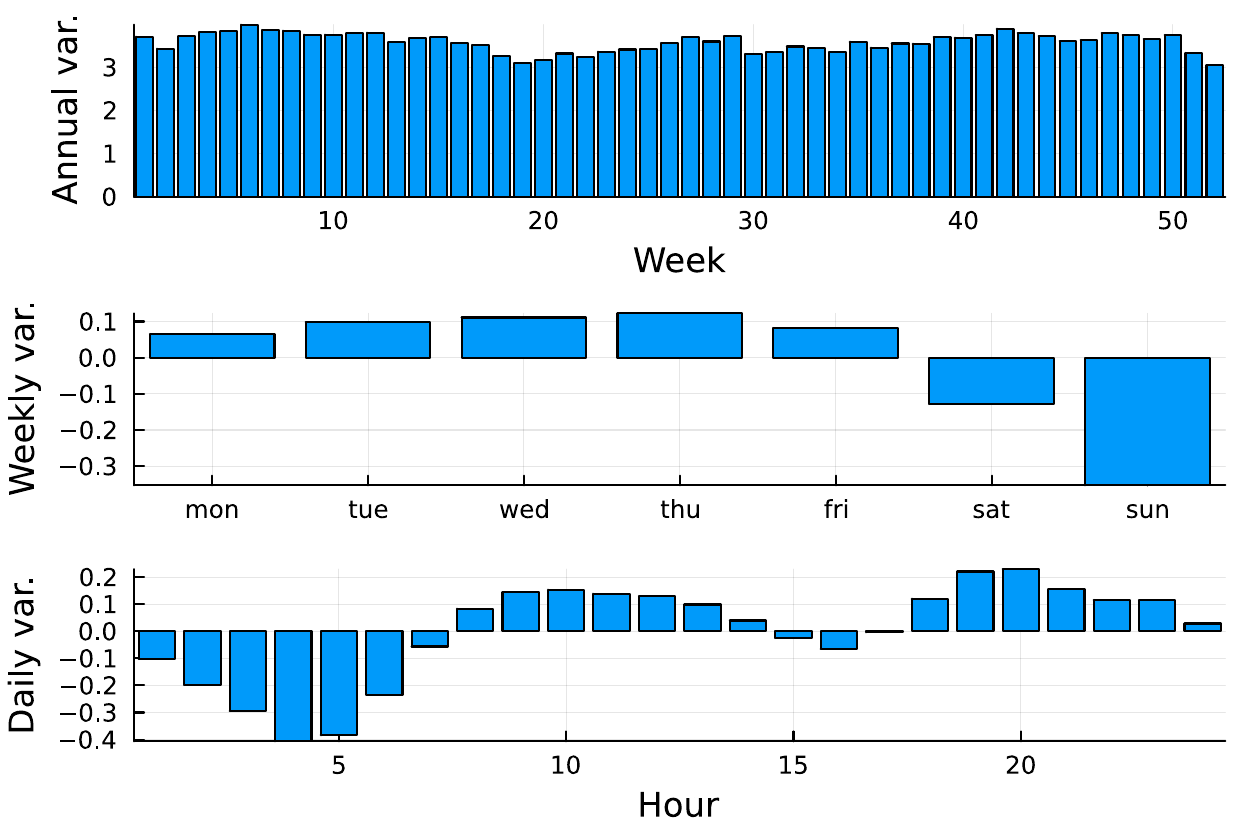}
  \caption{Annual, weekly and daily variations. \label{fig:mean-reversion}}
\end{figure}

\begin{figure}[!ht]
  \centering
  \begin{tabular}{cc}
    \includegraphics[width=.55\textwidth]{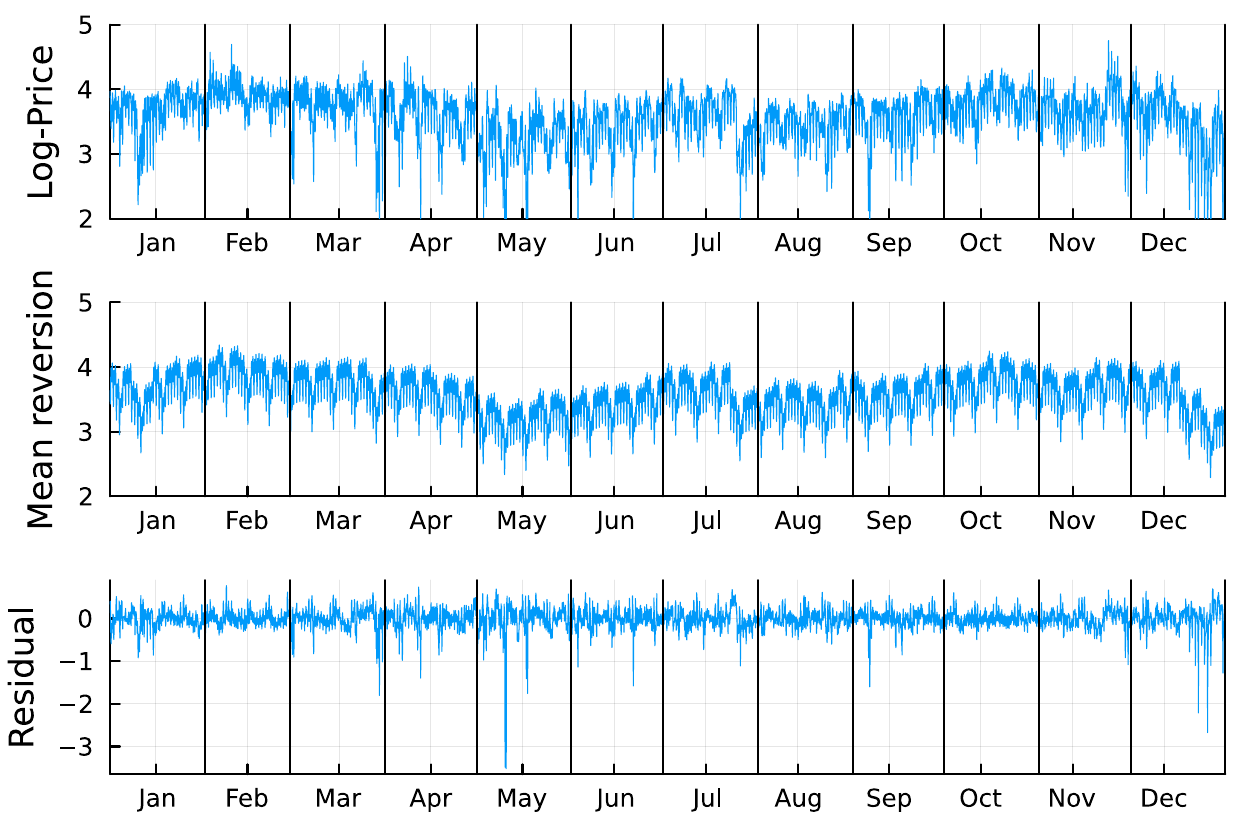}
    \includegraphics[width=.35\textwidth]{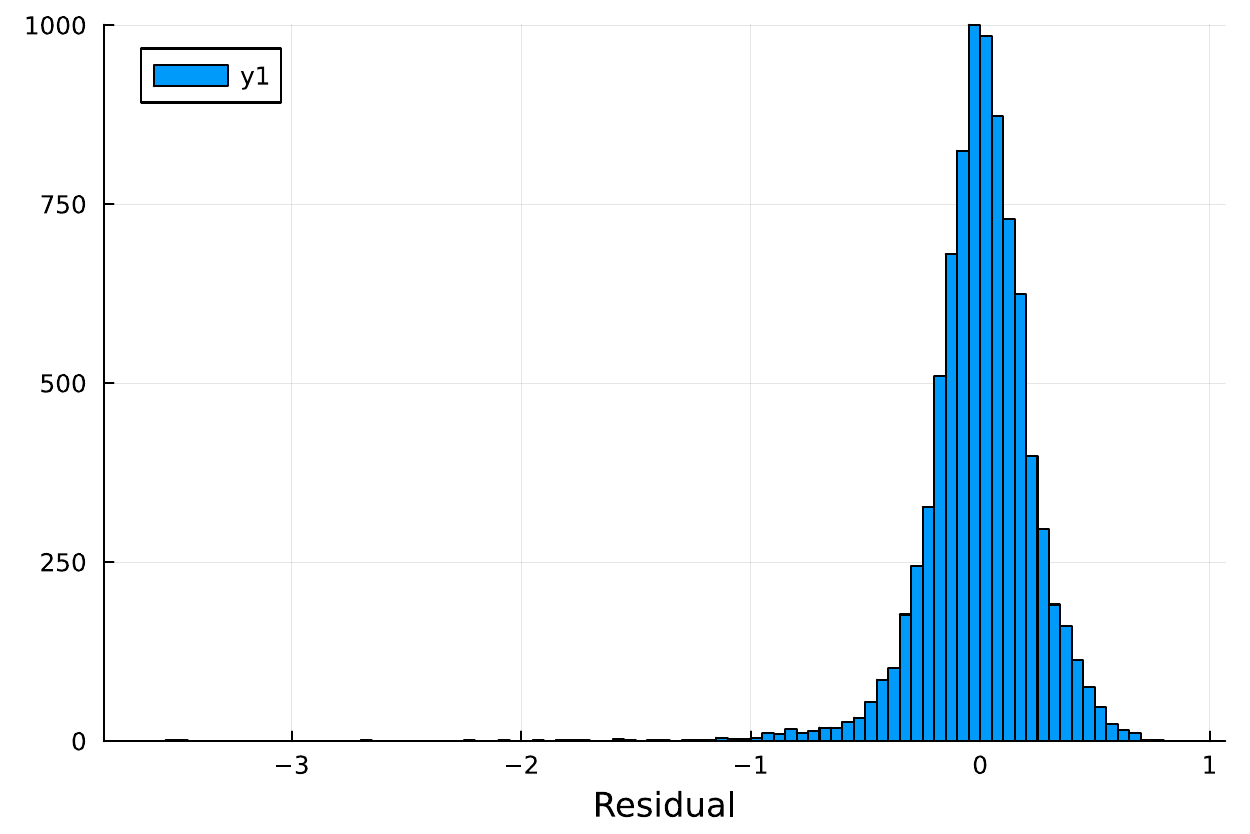}
  \end{tabular}
  \caption{Left: logarithmic spot price, average and residuals. Right: histogram of the residual of the time series model.
    \label{fig:epex-price}}
\end{figure}

The top plot in Figure~\ref{fig:epex-price} displays hourly logarithmic
electricity prices obtained from European Power Exchange (EPEX) during the year
of 2016. As a first step, we calculate the average cyclical variations of
electricity prices on the annual, weekly and daily levels by doing the
following:

\begin{enumerate}
\item We compute the average of the logarithmic price over each week of the year
  and then subtract the averages from the hourly values to remove the annual
  variation from the data;
\item Using the modified data from Item~1, we compute the average for each day
  of the week and subtract these from the data to remove the weekly variation
  from the data;
\item Using the modified data from Item~2, we compute the average for each hour
  of the day and subtract this from the data to remove the daily variation.
\end{enumerate}
Figure~\ref{fig:mean-reversion} plots the resulting average annual, weekly, and
daily patterns. Adding the three averages together gives an overall average
evolution of the logarithmic electricity price. Figure~\ref{fig:epex-price}
plots the resulting average as well the residual obtained by subtracting the
average from the original log-price data.

We denote the average logarithmic spot price by $m_t$ and model the deviations
\[
  \xi_t := \ln s_t-m_t
\]
by the mean-reverting Ornstein-Uhlenbeck process
\begin{equation}\label{eq:ou}
  \Delta\xi_{t+1} = -a\xi_t + \sigma\epsilon_t,
\end{equation}
where the constant $a$ and $\sigma$ are estimated from the data and
$(\epsilon_t)_{t=0}^T$ are independent Gaussians with zero mean and unit variance. This
results in a Markovian spot price model that will be used in the numerical
illustrations in Sect~\ref{sec:comp}.  Again, more sophisticated models could be
fitted and employed. 

\section{Computational results}\label{sec:comp}

In Section~\ref{subsec:Data}, we provide numerical data for the model presented
in Section~\ref{sec:example}.  In Section~\ref{sec:conv}, we analyze the
convergence of the SDDP algorithm.  Section~\ref{Sensitivity_analysis} is
devoted to sensitivity analysis of the numerical results with respect to
different parameters.

\subsection{Data}
\label{subsec:Data}

\begin{figure}[!ht]
  \centering
  \includegraphics[width=.7\textwidth]{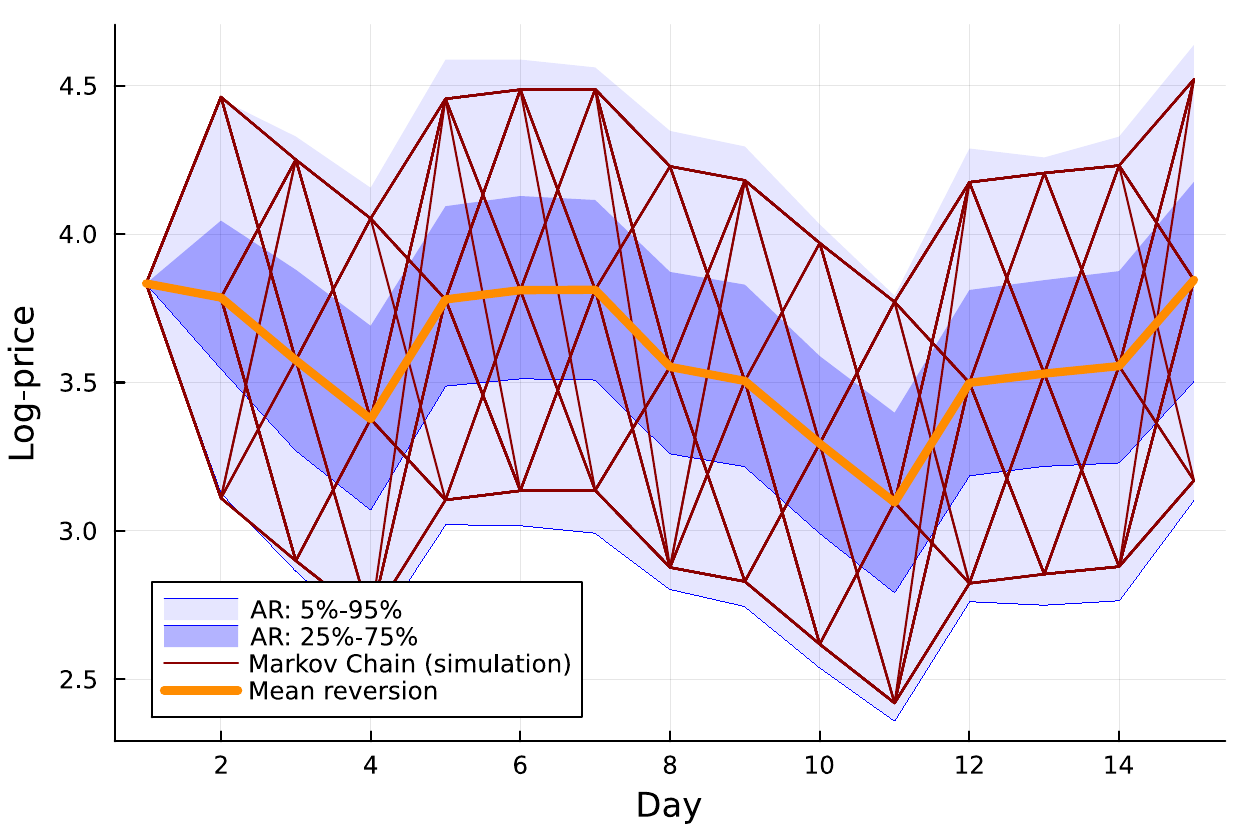}
  \caption{AR-1 process and its discretization as a Markov chain with 3 quadrature points per period \label{fig:markov-chain}}
\end{figure}

In the numerical studies below, we consider the problem of storage management
over one month with daily decisions so $T=30$. We discretize the one-dimensional
Markov process $\xi$ from Section~\ref{sec:pricemodel} by the procedure described
in Section~\ref{subsec:dmp}. We choose the sampling densities $\phi_t$ to be
Gaussian with zero mean and variance $\sigma^2$. The corresponding measures are then
approximated by Gauss-Hermite quadrature. For simplicity, we use the same
number~$N$ of quadrature points in each stage $t$. Figure~\ref{fig:markov-chain}
illustrates the discretization of the log-price over~15 days with $N=3$
quadrature points per time.

We solve the resulting discretized Bellman equations with the Julia package {\tt
  SDDP.jl}~\cite{dowson2021sddp} which implements the Markov chain SDDP
algorithm outlined in Section~\ref{subsec:sddp}.


\subsection{Convergence of the SDDP algorithm}\label{sec:conv}
We study the performance of SDDP as we increase the number $N$ of quadrature
points per period. For each $N$, we run 500 iterations of SDDP. The risk
aversion parameter $\rho$ is set to $10^{-4}$. The respective computation times are
given in Table~\ref{tab:computation_times}; they depend roughly linearly on~$N$.

\begin{table}[!ht]
  \centering
  \begin{tabular}{lrrrrrrr}
    \hline
    $N$ & 1 & 2 & 4 & 8 & 16 & 32 & 64 \\
    \hline
    Exec time (s) & 4.1 & 10.6 & 21.6 & 42.6 & 83.0 & 144.0 & 251.0  \\
    \hline
  \end{tabular}
  \caption{Time to compute 500 iterations
    of SDDP as we increase the number~$N$ of quadrature points in the discretization of the AR process. All times are given
    in seconds \label{tab:computation_times}}
\end{table}

\begin{figure}[!ht]
  \centering
  \begin{tabular}{cc}
    \includegraphics[width=.45\textwidth]{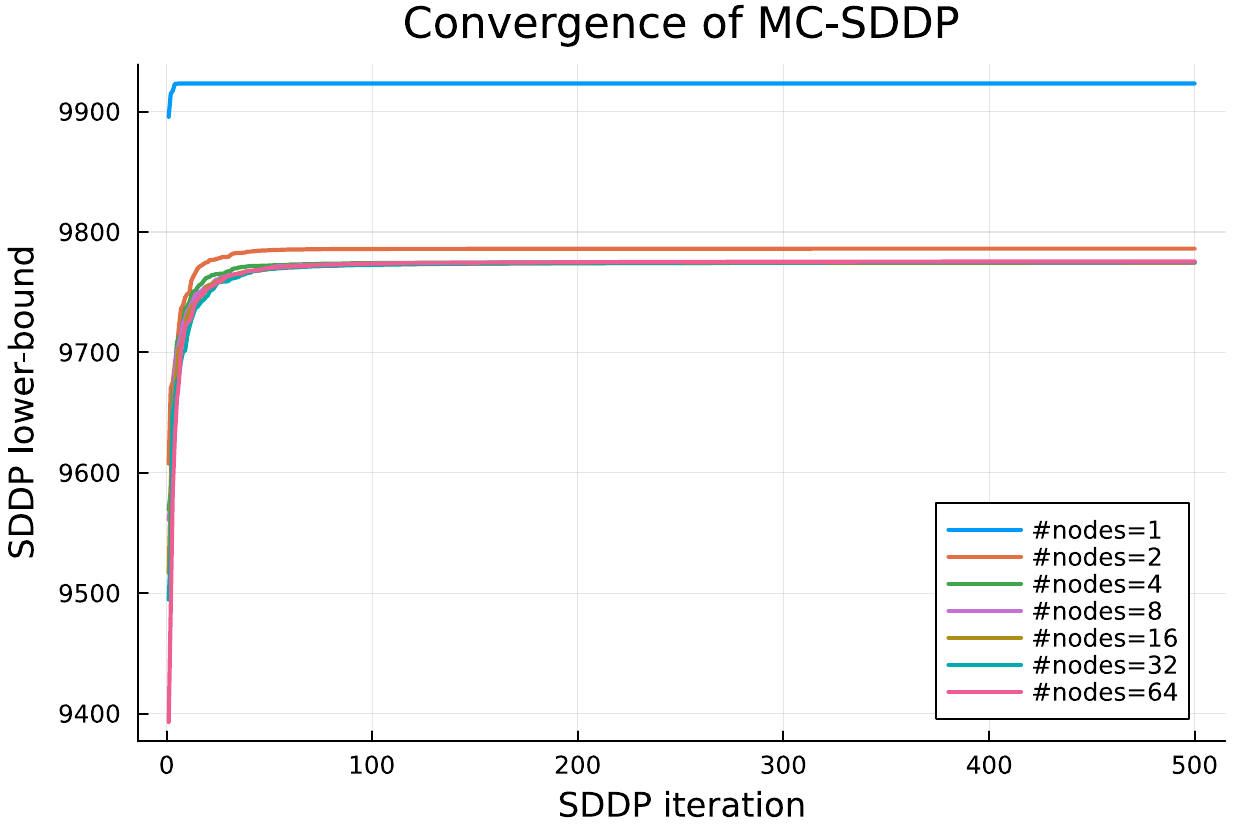}
    &
      \includegraphics[width=.45\textwidth]{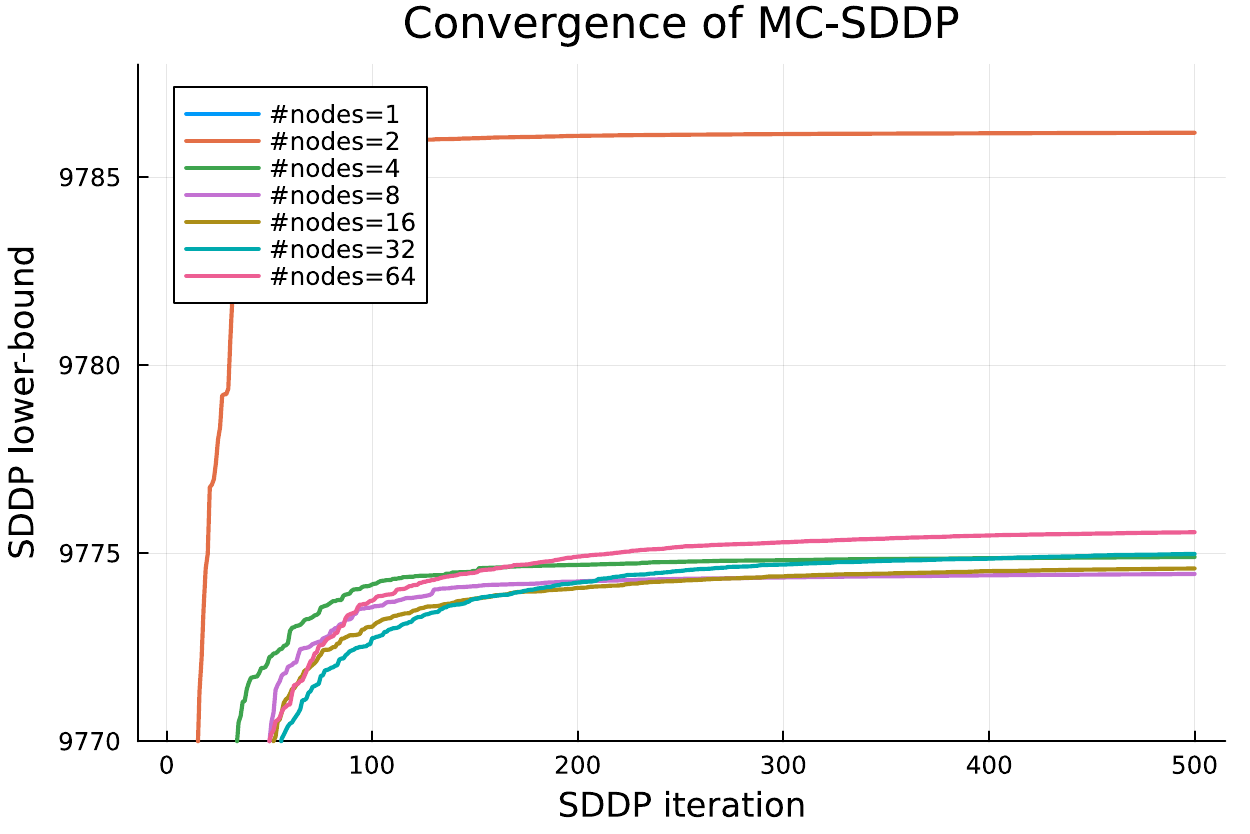}
  \end{tabular}
  \caption{Convergence of SDDP as we increase the number of quadrature points per period. The plot on the right is a zoomed version of the plot on the left.
  \label{fig:conv-sddp}}
\end{figure}

The convergence of the lower bounds is depicted in
Figure~\ref{fig:conv-sddp}. In the deterministic case where $N=1$, the algorithm
converges in less than 10 iterations. For $N\ge 2$, the convergence gets slower
but, interestingly, for $N=64$, the lower bound increases quicker than for
smaller values of $N$.

In order to evaluate the quality of the cost-to-go functions provided by the
SDDP algorithm, we perform an out-of-sample evaluation along 1,000 scenarios
randomly drawn from the continuously distributed Markov process described in
Section~\ref{sec:pricemodel}. In each scenario, we implement the control
obtained as in the forward pass of the SDDP algorithm as described in
Section~\ref{subsec:sddp}. However, since the random out-of-sample values of
$\xi_t$ will not coincide with the Markov states $s_t^i$ in the finite-state
approximation, we use the control obtained in the state $s_t^i$ nearest to the
sampled value of $\xi_t$.

\begin{figure}[!ht]
  \centering
  \includegraphics[width=.7\textwidth]{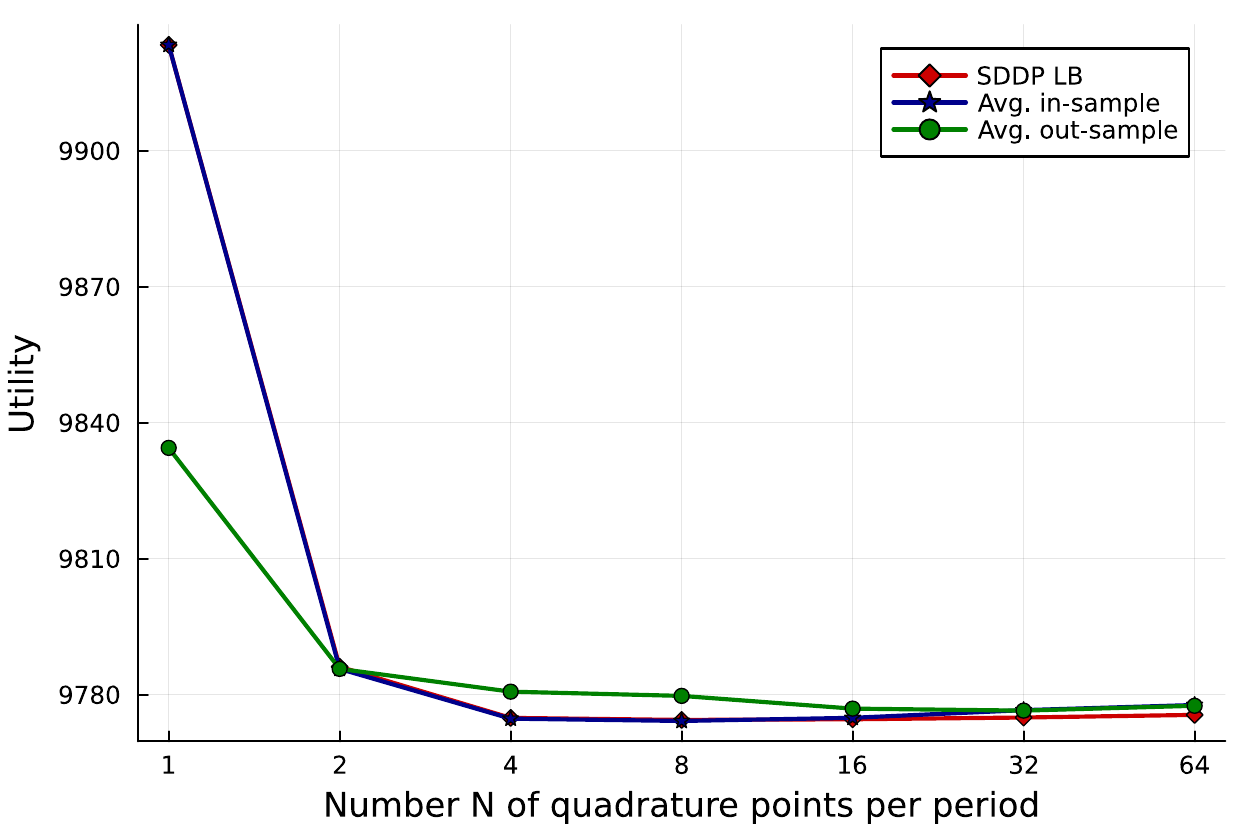}
  \caption{
    Optimum value estimates for increasing number of quadrature points. The red and blue values are given by SDDP while the green
    values are obtained with out of sample simulation using the optimized Bellman functions.
    \label{fig:assess-lb}}
\end{figure}

Figure~\ref{fig:assess-lb} plots the out-of-sample estimates of the objective
values together with the lower bounds found after 500 iterations of the SDDP
algorithm. The figure also plots ``in-sample" averages obtained by evaluating
the objective along scenarios randomly drawn from the finite-state Markov chain
approximation of the continuously distributed Markov process $\xi$.  The
out-of-sample averages decrease and tend towards the SDDP lower bound as we
increase $N$. This is expected, as the Markov discretization is an approximation
of the reference AR-1 process.

\begin{figure}[!ht]
  \centering
  \includegraphics[width=.9\textwidth]{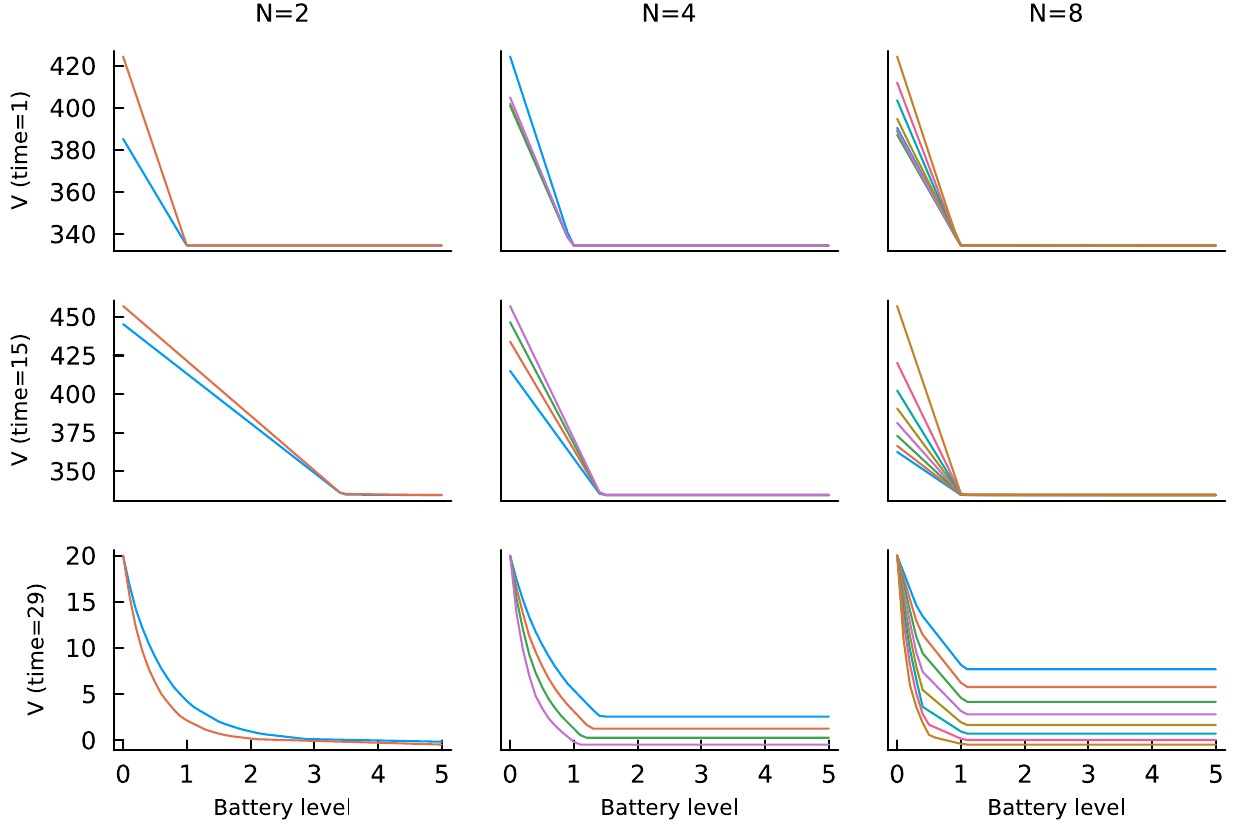}
  \caption{SDDP value functions as we increase the number~$N$ of quadrature points per period. The different rows correspond to different decision stages in the problem.
    \label{fig:value-function}}
\end{figure}

Figure~\ref{fig:value-function} plots some of the cost-to-go functions obtained
after 500 iterations of the SDDP algorithm with the number $N$ of quadrature
points equal to 2, 4 and 8. The figure plots the cost-to-go functions at three
different times: $t=1$, $t=15$ and $t=30$.


\subsection{Sensitivity analysis}
\label{Sensitivity_analysis}
As observed in Section~\ref{sec:conv}, the SDDP algorithm gives fairly accurate
results already with $N = 8$ quadrature points per period. From now on, we fix
$N=8$ and proceed to study the sensitivity of the solutions with respect to the
parameters of the problem. We study both the solution of the optimization
problem~\eqref{sps} as well as the indifference prices obtained from using the
optimum value function of~\eqref{sps} as described in
Section~\ref{sec:idp}. Because, in the present example, the interest rates are
constant, the exponential utility allows us to find the indifference price
analytically after the solution of just one instance of the battery management
problem.

\begin{figure}[!ht]
  \centering
  \includegraphics[width=.75\textwidth]{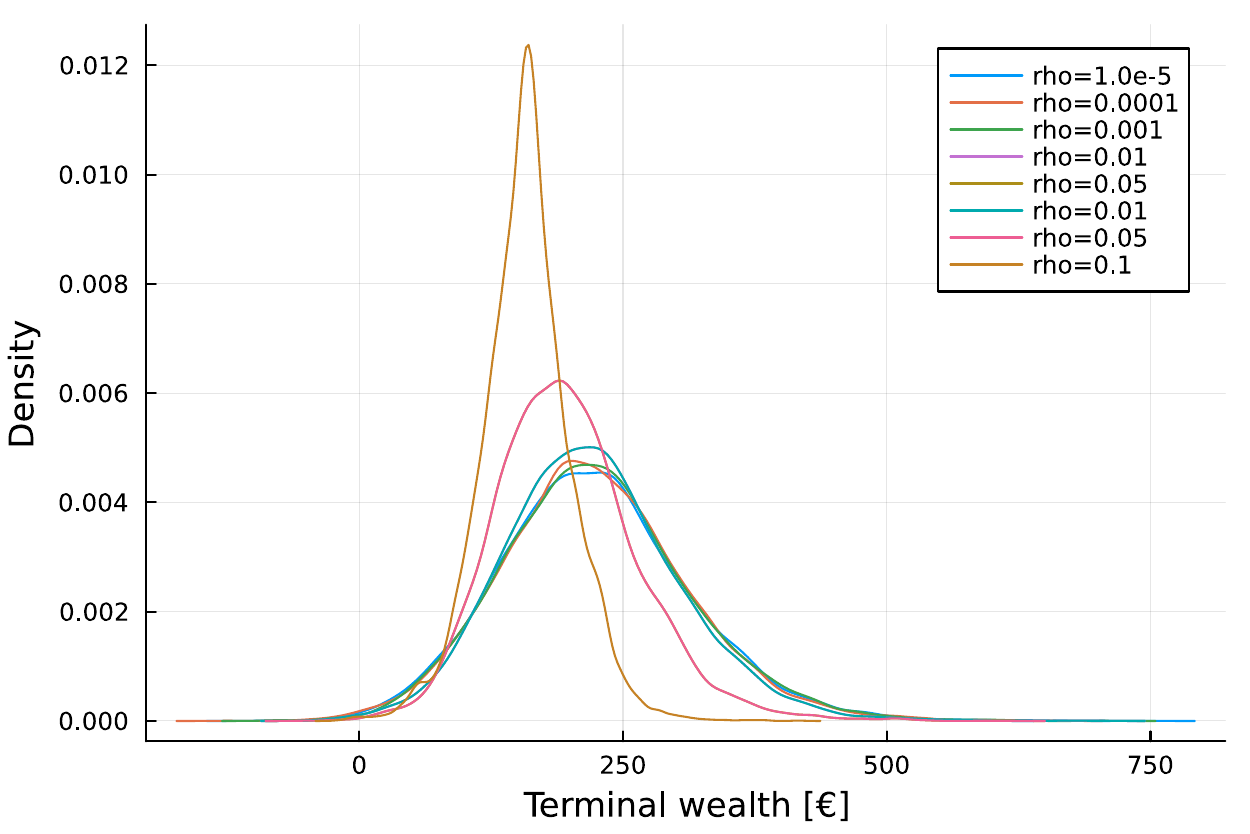}
  \caption{
    Kernel density plots of the terminal wealth in out-of-sample simulations using 10,000 scenarios for different
    values of the risk aversion parameter $\rho$ in the exponential utility function~\eqref{eq:exp}. \label{fig:sensitivity:risk}
  }
\end{figure}

Figure~\ref{fig:sensitivity:risk} displays kernel-density plots of the terminal
wealth obtained with strategies optimized for different risk aversion. As
expected, higher risk aversion results in thinner left tail of the terminal
wealth distribution.

\begin{figure}[!ht]
  \centering
  \includegraphics[width=.7\textwidth]{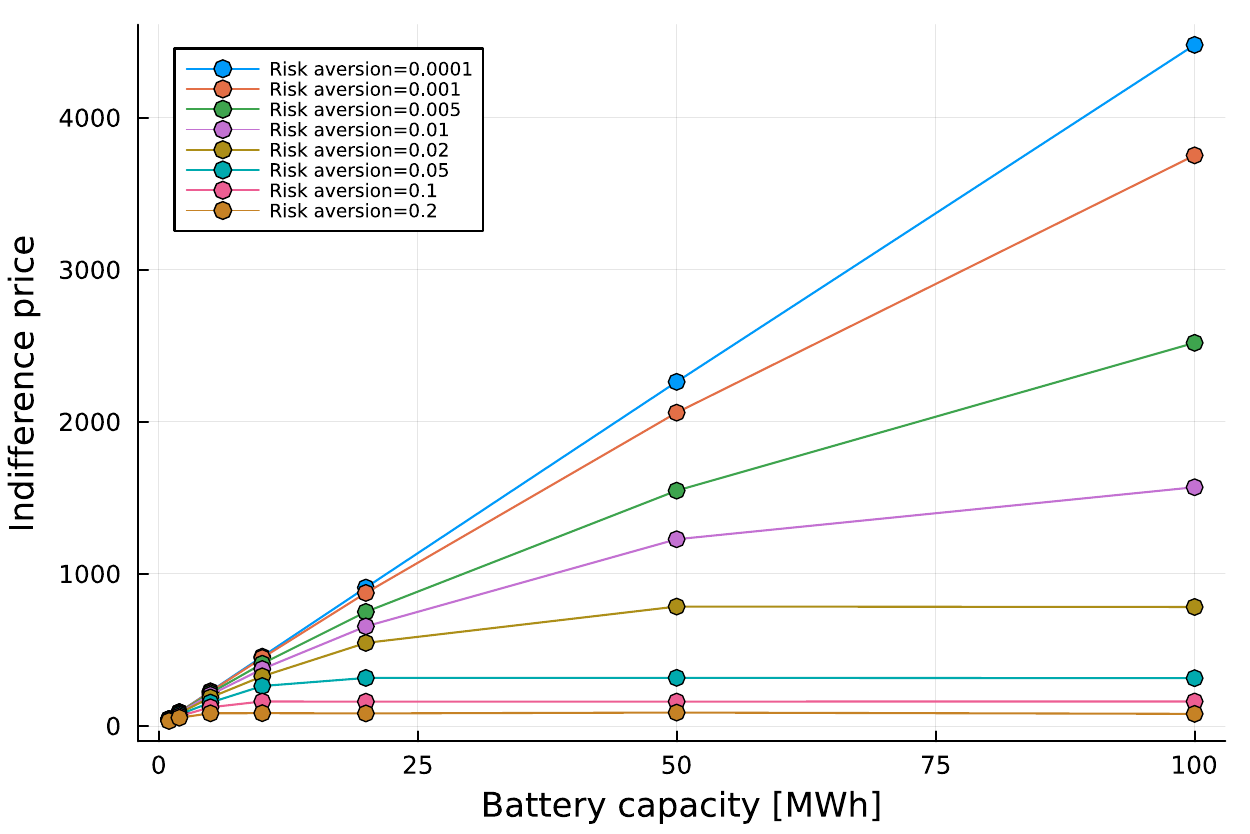}
  \caption{Indifference prices (in euros) of storages with different storage capacities.
    Different colors correspond to agents with different risk aversion parameters.
    Prices are given in euros and capacity in MWh.
  \label{fig:indifference:capacity}}
\end{figure}
\begin{figure}[!ht]
  \centering
  \includegraphics[width=.7\textwidth]{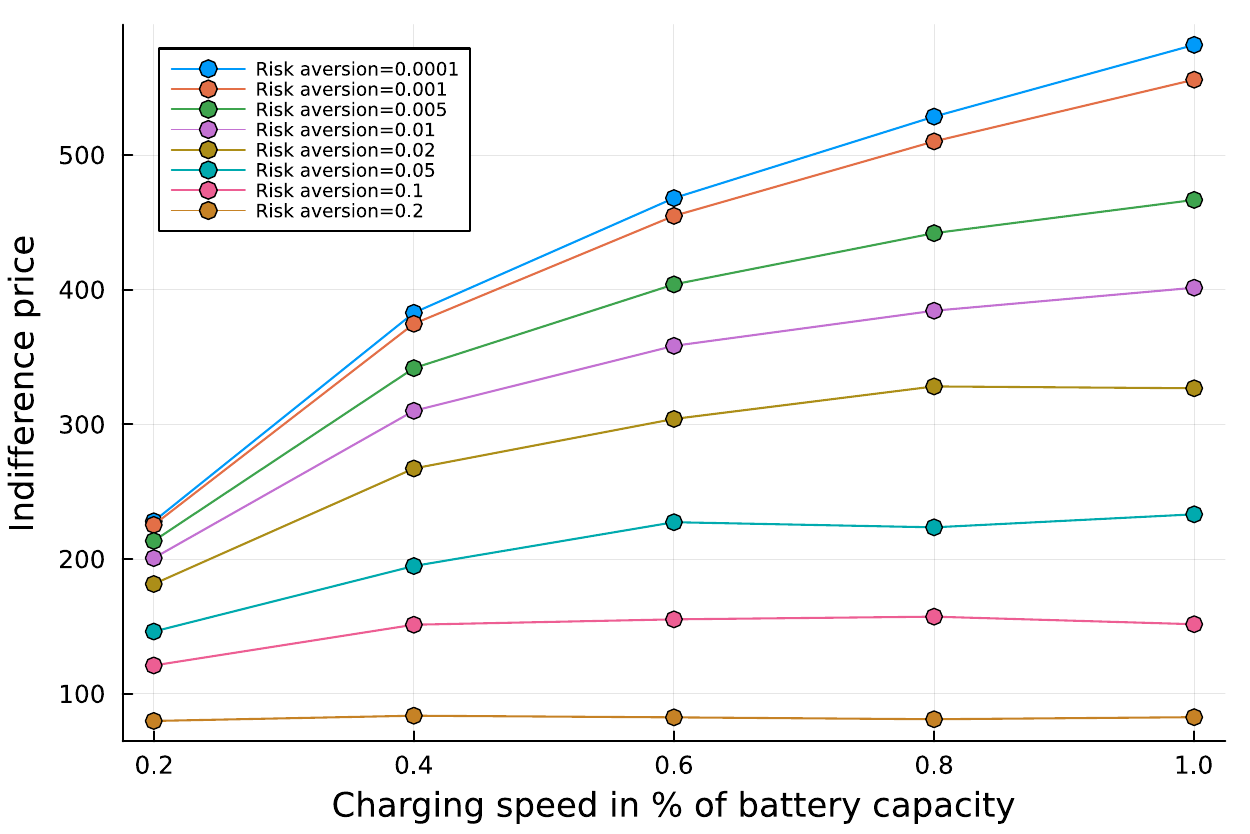}
  \caption{Indifference prices (in euros) of storages with different charging speeds.
    The horizontal axis gives the maximum daily change in stored energy in units of the total storage capacity.
    In particular, charging speed of $1$ means that the charging speed is unconstrained.
    Different colors correspond to agents with different risk aversion parameters.
  \label{fig:indifference:rate}}
\end{figure}

We will finish this section by studying sensitivities of the indifference price
with respect to the parameters of the
model. Figure~\ref{fig:indifference:capacity} plots indifference prices for
storages with different capacities. As one might expect, the valuation of a
storage increases with the storage capacity. With higher risk aversion, however,
the valuations seem to saturate for larger capacities. This is quite natural as
a highly risk-averse agent will avoid large positions in stored electricity due
to the price risk. Figure~\ref{fig:indifference:rate} plots the indifference
prices for storages with different charging speeds. Again, the valuations are
more sensitive for less risk averse agents. Both figures also give the prices
for agents with varying risk aversion parameters. Finally,
Figure~\ref{fig:indifference:rate} gives indifference prices for varying values
of the volatility parameter $\sigma$ in the Ornstein-Uhlenbeck
process~\eqref{eq:ou}. The standard deviation can be thought of as the agent's
uncertainty concerning the evolution of electricity prices. For all levels of
risk aversion, higher volatility makes the storage more valuable to the agent.

\begin{figure}[!ht]
  \centering
  \includegraphics[width=.7\textwidth]{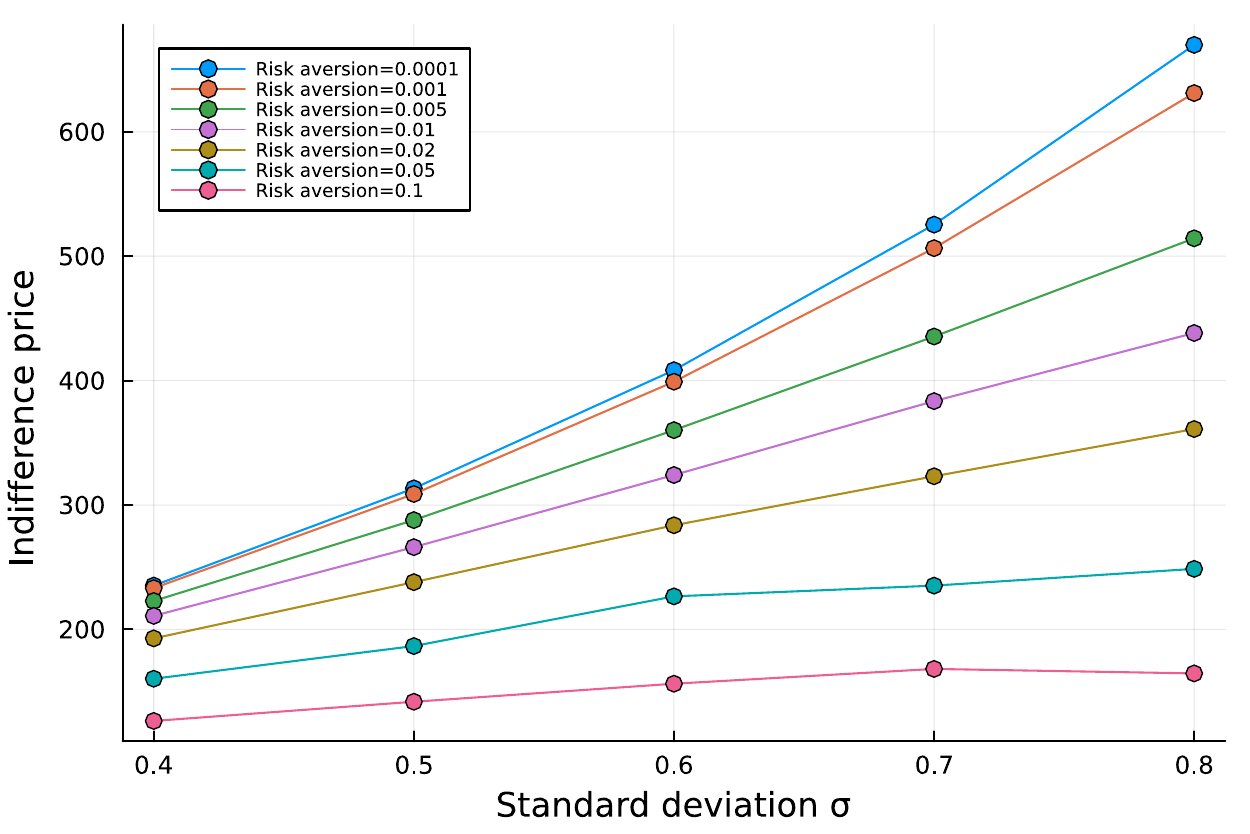}
  \caption{Indifference prices (in euros) obtained with different values of the volatility parameter $\sigma$
    in the Ornstein-Uhlenbeck process~\eqref{eq:ou}.
    Different colors correspond to different risk aversions.
  \label{fig:indifference:std}}
\end{figure}

\bibliographystyle{plain}
\bibliography{sp}

\end{document}